\documentclass[regno]{amsart}
\usepackage{amsmath,amssymb,amscd}

\newtheorem{Thm}{Theorem}[section]

\newtheorem{Rem}{Remark}[section]

\newtheorem{Example}{Example}[section]

\newtheorem{Conjecture}{Conjecture}

\usepackage{comment}

\numberwithin{equation}{section}

\begin{document}
\title{On the generic Conley conjecture}
\author{Yoshihiro Sugimoto}
\date{}
\maketitle

\begin{abstract}
In this paper, we treat an open problem related to the number of periodic orbits of Hamiltonian diffeomorphisms on closed symplectic manifolds, so-called (generic) Conley conjecture. Generic Conley conjecture states that generically Hamiltonian diffeomorphisms have infinitely many simple contractible periodic orbits. We prove generic Conley conjecture for very wide classes of symplectic manifolds.
\end{abstract}

\section{Introduction and main results}
In this section, we briefly explain the main theme of this paper. The precise definitions and notations are given in the next section. The information of periodic orbits of Hamiltonian diffeomorphisms is very important in Hamiltonian dynamics. Conley conjecture was originally stated for Hamiltonian diffeomorphisms on the standard torus ${(\mathbb{T}^{2n}, \omega_{0})}$ (\cite{C}). It states that any Hamiltonian diffeomorphism on ${(\mathbb{T}^{2n}, \omega_{0})}$ has infinitely many simple contractible periodic orbits (Simple means it is not iterated periodic orbit of lower period.). It is easy to see that this conjecture can not be generalized to any closed symplectic manifolds. For example, an irrational rotation on the standard sphere ${S^2\subset \mathbb{R}^3}$ has only two contractible periodic orbits, the north pole and the south pole.

However, Conley conjecture was proved for wide classes of closed symplectic manifolds. For example, Conley conjecture holds on symplectically aspherical manifolds, negatively monotone symplectic manifolds and symplectic manifolds with vanishing spherical Chern class(\cite{G1,GG2, GG3,GG4,H, SaZ}). So, today's Conley conjecture is a conjecture that every Hamiltonian diffeomorphism has infinitely many simple contractible periodic orbits on "almost all" closed symplectic manifolds.

Another variant of above Conley conjecture is so-called generic Conley conjecture (\cite{GG,GG3,GG4}). Generic conley conjecture states that "almost all" Hamiltonian diffeomorphisms have infinitely many simple contractible periodic orbits on every closed symplectic manifold. Conley conjecture and generic Conley conjecture state that Hamiltonian diffeomorphims with finitely many simple periodic orbits (like the irrational rotation on the sphere ${S^2}$) are very rare. In summary, we have the following two conjectures.

\begin{Conjecture}[(generic) Conley conjecture]\label{ConjectureA}
\begin{enumerate}
\item On "almost all" closed symplectic manifolds, every Hamiltonian diffeomorphism has infinitely many simple contractible periodic orbits.
\item On every closed symplectic manifolds, almost all Hamiltonian diffeomorphisms have infinitely many simple contractible periodic orbits.
\end{enumerate}
\end{Conjecture}

 In this paper, we study Conjecture \ref{ConjectureA} ${(2)}$, the generic Conley conjecture. The statement of our main result is stated as follows. 

\begin{Thm}
Let ${(M,\omega)}$ be a ${2n}$-dimensional closed symplectic manifold and let ${N\in \mathbb{N}\cup \{\infty\}}$ be the minimum Chern number of ${(M,\omega)}$. Assume that ${(M,\omega)}$ satisfies at least one of the following conditions.
\begin{enumerate}
\item ${n}$ is odd.
\item ${H_{odd}(M:\mathbb{Q})\neq 0}$
\item ${N>1}$
\end{enumerate}
Then, there is a ${C^{\infty}}$-dense and ${C^{\infty}}$-residual ($=$contains a countable intersection of ${C^{\infty}}$-open dense subsets) subset ${\mathcal{U}\subset \textrm{Ham}(M,\omega)}$ such that any element of ${\mathcal{U}}$ has infinitely many simple contractible periodic orbits.
\end{Thm}

Note that the above conditions ${(1)}$, ${(2)}$ and ${(3)}$ cover almost all closed symplectic manifolds.

\begin{Rem}
The case (2) of Theorem 1.1 was also proved in Proposition 1.6 in \cite{GG}. The case (3) is a generalization of Theorem 1.2 in \cite{GG} where Ginzburg and G\"{u}rel proved generic Conley conjecture for ${N\ge n+1}$. The proof of Proposition 1.6 in \cite{GG} was an application of Birkhoff-Moser fixed point theorem and the proof of Theorem 1.2 in \cite{GG} was an application of "resonance relation" proved in \cite{GK}. Our proof of Theorem 1.1 is a modification of the former proof.
\end{Rem}

\section*{Acknowledgement}
This work was carried out during my stay as a research fellow in National Center for Theoretical Sciences. The author thanks NCTS for a great research atmosphere and many supports. He also gratefully acknowledges his teacher Kaoru Ono for continuous supports and Victor L. Ginzburg for checking the draft and giving me comments and advices.

\section{Preliminaries}
In this section, we explain notations and terminologies used in this paper.

\subsection{Elementary notations}
Let ${(M,\omega)}$ be a symplectic manifold, so $M$ is a finite dimensional $C^{\infty}$-manifold and ${\omega \in \Omega^2(M)}$ is a symplectic form on $M$. In this paper, we always assume that $M$ is a closed manifold.

For any $C^{\infty}$ function $H\in C^{\infty}(M)$, we define the Hamiltonian vector field ${X_H}$ by the following relation.

\begin{equation*}
\omega(X_H, \cdot)=-dH
\end{equation*}

We can also consider ${S^1}$-dependent ($=$$1$-periodic) Hamiltonian function $H$ and Hamiltonian vector field $X_H$ by the same formula. The time $1$  flow of $X_H$ is called a Hamiltonian diffeomorphism generated by $H$. We denote this flow by $\phi_H$. The set of all Hamiltonian diffeoomorphisms is called Hamiltonian diffeomorphism group and we denote the Hamiltonian diffeomorphism group of ${(M,\omega)}$ by ${\textrm{Ham}(M,\omega)}$.

\begin{equation*}
\textrm{Ham}(M,\omega)=\{\phi_H \ | \ H\in C^{\infty}(S^1\times M)\}
\end{equation*}

We also consider "iterations" of $H$ and ${\phi_H}$. For any integer ${k\in \mathbb{N}}$, we define ${H^{(k)}}$ as follows.

\begin{equation*}
H^{(k)}=kH(kt, x)
\end{equation*}
It is straightforward to see that ${\phi_{H^{(k)}}=(\phi_H)^k}$. Let ${P^l(H)}$ be the space of $l$-periodic contractible periodic orbits of $X_H$. 

\begin{gather*}
P^l(H)=\{x:S_l^1\rightarrow M \ | \ \dot{x}(t)=X_{H_t}(x(t)), x :\textrm{contractible}  \}  \\
S_l^1=\mathbb{R}/l\cdot \mathbb{Z}
\end{gather*}
It is also straightforward to see that there is one to one correspondence between ${P^k(H)}$ and ${P^1(H^{(k)})}$. We abbreviate ${P^1(H)}$ to ${P(H)}$. A $l$-periodic orbits ${x\in P^l(H)}$ is called simple if there is no $l'$-periodic orbits ${y\in P^{l'}(H)}$ which satisfies the following conditions.

\begin{gather*}
l=l'\cdot m  \ \ \ \ (l', m\in \mathbb{N}) \\ 
x(t)=y(\pi_{l,l'}(t))
\end{gather*} 
Here ${\pi_{l,l'}:S_l\rightarrow S_{l'}}$ is the natural projection. So a periodic orbit is simple if and only if it is not iterated periodic orbit of lower period.

Next, we explain the definition of the minimum Chern number $N$. A symplectic manifold ${(M,\omega)}$ becomes an almost complex manifold, and its tangent bundle has a natural first Chern class ${c_1(TM)\in H^2(M:\mathbb{Z})}$. The minimum Chern number $N\in \mathbb{N}\cup \{+\infty\}$ is the positive generator of ${c_1(TM)|_{\pi_2(M)}}$. Note that if the image is zero, ${N}$ is defined by ${N=+\infty}$.

\subsection{Floer homology and degrees of periodic orbits}
In this subsectin, we explain basic notations of Floer homology theory and Conley-Zehnder index of periodic orbit. Let $H$ be a $1$-periodic Hamiltonian function. We call $H$ is non-degenerate if the differential map ${d\phi_H:TM_x\rightarrow TM_x}$ does not has $1$ as an eigenvalue for any fixed point ${x\in \textrm{Fix}(\phi_H)}$. Note that $H$ is non-degenerate if and only if  ${\textrm{graph}(\phi_H)\subset M\times M}$ is transverse to the diagonal ${\Delta_M \subset M\times M}$.

We construct Novikov covering of ${P(H)}$ as follows.

\begin{equation*}
\widetilde{P(H)}=\{(u,x) \ | \ x\in P(H), u:D^2\rightarrow M, \partial u=x \}/\sim
\end{equation*}
where $D^2$ is the two dimensional disc ${D^2\subset \mathbb{R}^2}$ and the equivalence relation ${\sim}$ is defined as follows.
\begin{gather*}
(u,x)\sim(v,y)\Longleftrightarrow \begin{cases} x=y  \\ \omega(u\sharp \overline{v})=0  \\ c_1(u\sharp \overline{v})=0
\end{cases}
\end{gather*}
Here ${\overline{v}}$ is the disc with the opposite orientation on the domain and ${u\sharp \overline{v}}$ is the glued sphere. Each ${[u,x]\in \widetilde{P(H)}}$ has a Conley-Zehnder index ${\mu_{CZ}([u,x])\in \mathbb{Z}}$. We normalize ${\mu_{CZ}}$ so that Conley-Zehnder index of a local maximum of a ${C^2}$-small Morse function is equal to $n$. Conley-Zehnder index gives a grading of Floer chain complex and Floer homology. We also have the action functional ${A_H}$ on ${\widetilde{P(H)}}$ as follows.

\begin{equation*}
A_H([u,x])=-\int_{D^2}u^*\omega+\int_0^1H(t, x(t))dt
\end{equation*}
Then Floer chain complex ${CF_*(H)}$ is defined as follows.

\begin{equation*}
CF_*(H)=\bigg\{\sum_{z\in \widetilde{P(H)}}a_z\cdot z \ \bigg| \ a_z\in \mathbb{Q}, \forall C\in \mathbb{R}, \sharp\{z\in \widetilde{P}(H) \ | \ a_z\neq 0, A_H(z)>C\}< \infty \bigg\}
\end{equation*}
The boundary operator $d_F$ has the following form.
\begin{equation*}
d_F(z)=\sum_{w\in \widetilde{P(H)}}n(z,w)w
\end{equation*}
The coefficient ${n(z,w)\in \mathbb{Q}}$ is the number of solutions of the following Floer equation modulo the natural ${\mathbb{R}}$-action (\cite{FO2,LT}). Let ${J_t}$ be an almost complex structure on $M$ parametrized by ${t\in S^1}$. 

\begin{gather*} 
z=[v_-,x_-], w=[v_+,x_+]  \\
u:\mathbb{R}\times S^1\longrightarrow M  \\
\partial _su(s,t)+J_t(u(s,t))(\partial _tu(s,t)-X_{H_t}(u(s,t)))=0  \\
\lim_{s\to -\infty}u(s,t)=x_-(t), \lim_{s\to +\infty}u(s,t)=x_+(t), (v_-\sharp u, x_+)\sim (v_+,x_+)
\end{gather*}

Floer homology ${HF_*(H)}$ is the homology of the chain complex ${(CF_*(H),d_F)}$. We introduce the notion of Novikov ring of ${(M,\omega)}$. We define an abelian groug $\Gamma$ by

\begin{equation*}
\Gamma=\frac{\pi_2(M)}{\textrm{Ker}\omega \cap \textrm{Ker}c_1}
\end{equation*}
where ${\omega:\pi_2(M)\rightarrow \mathbb{R}}$ is the integration of the symplectic form $\omega$ and ${c_1:\pi_2(M)\rightarrow \mathbb{Z}}$ is the integration of the first Chern class. We define the degree of ${u\in \Gamma}$ by ${-2c_1(u)}$. Novikov ring ${\Lambda_{(M,\omega)}}$ is defined by the set of possibly infinite sums of $\Gamma$ with suitable convergence as follows.

\begin{gather*}
\Lambda_{(M,\omega)}=\bigg\{\sum_{u\in \Gamma}a_u\cdot u \ \bigg| \ a_u\in \mathbb{Q},\forall C\in \mathbb{R},\sharp\{u\in \Gamma \ | \ a_u\neq 0, \omega(u)<C\}<\infty  \bigg\}
\end{gather*}

Then Floer homology is isomorphic to the singular homology group with Novikov ring coefficient (\cite{FO2,LT}).

\begin{equation*}
HF_*(H)\cong H_{*-n}(M:\mathbb{Q})\otimes \Lambda_{(M,\omega)}
\end{equation*}

\section{Generic Conley conjecture}

We prove Theorem 1.1 in this section. Throughout this section, we assume that ${(M,\omega)}$ is a ${2n}$-dimensional closed symplectic manifold with minimum Chern number $N$ and it also satisfies at least one of the following conditions.

\begin{enumerate}
\item ${n}$ is odd.
\item ${H_{odd}(M:\mathbb{Q})\neq 0}$
\item ${N>1}$
\end{enumerate}

The purpose of this section is to construct a subset ${X\subset \textrm{Ham}(M,\omega)}$ and a family of subsets ${\{Y_k\subset \textrm{Ham}(M,\omega)\}}$ (${1\le k < +\infty}$) which satisfy the following conditions.

\begin{itemize}
\item ${X\subset \textrm{Ham}(M,\omega)}$ is a $C^{\infty}$-dense subset.
\item ${Y_k\subset \textrm{Ham}(M,\omega)}$ are $C^{\infty}$-open dense subsets.
\item Any element of $X$ has infinitely many simple contractible periodic orbits.
\item $X=\bigcap_{k=1}^{\infty}Y_k$ holds. 
\end{itemize}
The above conditions imply that $X$ is a ${C^{\infty}}$-residual subset of ${\textrm{Ham}(M,\omega)}$ and generically Hamiltonian diffeomorphisms have infinitely many simple contractible periodic orbits.

As in \cite{GG}, our proof is based on applications of Birkhoff-Moser fixed point theorem (local theory) and Floer homology theory (global theory). Roughly speaking, Birkhoff-Moser fixed point theorem guarantees infinitely many periodic orbits of a symplectic map near non-hyperbolic fixed points which satisfies some generic conditions. For the reader's convenience, we briefly recall the statement and properties of Birkhoff-Moser fixed point theorem.

\begin{Thm}[Birkhoff-Moser fixed point theorem \cite{Mo}]
Let $\phi$ be s symplectic map defined in an open neighborhood of the origin ($=p$) in ${(\mathbb{R}^{2n},\omega_0)}$ and the origin is a fixed point of ${\phi}$. Here $\omega_0$ is the standard symplectic form ${\sum_{i=1}^n x_i\wedge y_i}$ on ${\mathbb{R}^{2n}}$. Let ${\lambda_1,\cdots,\lambda_m,\lambda_1^{-1},\cdots,\lambda_m^{-1}}$ be the all eigenvalues of the differential map
\begin{equation*}
d\phi_p:T_pM\longrightarrow T_pM
\end{equation*}
on the unit circle in ${\mathbb{C}}$. Assume that $\phi$ satisfies the following conditions.
\begin{enumerate}
\item $m\ge 1$
\item $\prod_{k=1}^m\lambda_k^{j_k}\neq 1$ for ${1\le \sum_{k=1}^m|j_k|\le 4}$ 
\item The Taylor coefficient of $\phi$ up to order $3$ satisfies a non-degenerate condition.
\end{enumerate}
Then $\phi$ possesses infinitely many periodic orbits in any neighborhood of $p$.
\end{Thm}

The meaning of "non-degenerate" in ${(3)}$ is difficult to state briefly because its meaning becomes clear in the proof of the theorem. We just introduce an example of "non-degenerate" condition.

\begin{Example}[non-degeneracy condition \cite{Mo}]
Let ${\phi}$ be a symplectic map defined in a open neighborhood of the origin in ${(\mathbb{R}^{2n},\omega_0)}$ and the origin is a fixed point of ${\phi}$. Assume that $\phi$ can be written in the following form .
\begin{gather*}
\begin{cases}
\phi((x_1,\cdots,x_n,y_1,\cdots,y_n))=(x_1^{(1)},\cdots,x_n^{(1)},y_1^{(1)},\cdots,y_n^{(1)}) \\
x_k^{(1)}=x_k\cos \Phi_k-y_k\sin \Phi_k+f_k  \\
y_k^{(1)}=x_k\sin \Phi_k+y_k\cos \Phi_k+f_{k+n}  \\
\Phi_k =\alpha_k +\sum_{l=1}^n\beta_{kl}(x_l^2+y_l^2)
\end{cases}
\end{gather*}
The error terms ${f_k}$ are assumed to have vanishing derivatives up to order 3 at the origin. Then non-degeneracy means that the matrix ${(\beta_{kl})}$ is non-singular.
\end{Example}
Moser first proved Birkhoff-Moser fixed point theorem for the above special case. Then he proved that general cases can be reduced to this special case. So roughly speaking, "non-degenerate" means that it can be reduced to the above form so that the matrix ${(\beta_{kl})}$ is non-singular.

Birkhoff-Moser fixed point theorem (and its proof in \cite{Mo}) implies the following fact. Let ${x\in P(H)}$ be a non-degenerate contractible periodic orbit of a Hamiltonian function ${H\in C^{\infty}(S^1\times M)}$. We also assume that there is at least one eigenvalue of the differential map 

\begin{equation*}
d\phi_{H}:T_{x(0)}M\longrightarrow T_{x(0)}
\end{equation*}
on the unit circle and all eigenvalues on the unit circle are pairwise distinct.Then we can perturb $H$ to $\widetilde{H}$ near $x$ so that it satisfies the all required conditions in the statement of Birkhoff-Moser fixed point theorem. Moreover, these conditions are satisfied in a sufficiently small open neighborhood of ${\phi_{\widetilde{H}}}$ and hence all of them possesses infinitely many simple contractible periodic orbits.

We apply this observation to our proof of Theorem 1.1. Let ${\mathcal{H}_{sn}\subset \textrm{Ham}(M,\omega)}$ be the set of strongly non-degenerate Hamiltonian diffeomorphisms (strongly non-degenerate means any iteration of it is non-degenerate). We divide ${\mathcal{H}_{sn}}$ into the following three pairwise disjoint subsets.

\begin{gather*}
\mathcal{H}_{sn}^{(1)}=\Bigg\{\phi \in \mathcal{H}_{sn} \ \Bigg| \ \begin{matrix}\textrm{the\ number\ of simple\ contractible\ periodic\ orbits\ is\ finite } \\  \textrm{and\ all\ contractible\ periodic\ orbits\ are\ hyperbolic} \end{matrix}  \Bigg\}  
\end{gather*}
\begin{gather*}
\mathcal{H}_{sn}^{(2)}=\Bigg\{\phi \in \mathcal{H}_{sn} \ \Bigg| \ \begin{matrix}\textrm{the\ number\ of simple\ contractible\ periodic\ orbits\ is\ finite } \\  \textrm{and\ at\ least\ one\ contractible\ periodic\ orbit\ is\ non-hyperbolic} \end{matrix}  \Bigg\}
\end{gather*}
\begin{gather*}
\mathcal{H}_{sn}^{(3)}=\Bigg\{\phi \in \mathcal{H}_{sn} \ \Bigg| \ \phi \ \textrm{possesses infinitely many simple contractible periodic orbits}  \Bigg\}
\end{gather*}

First we prove that ${\mathcal{H}_{sn}^{(1)}}$ is empty. We fix ${\phi \in \mathcal{H}_{sn}^{(1)}}$ and let ${\{x_i,\cdots,x_l\}}$ be the set of all simple contractible periodic orbits of ${\phi}$ and let ${p_1,\cdots,p_l\in \mathbb{N}}$ be their periods. We also choose a common multiple ${k}$ of ${p_1,\cdots,p_l}$. Then all periodic orbits of ${\psi'=\phi^{(k)}}$ are $1$-periodic orbits and all of them are hyperbolic. For any capped periodic orbit ${\bar{z}}$, we have the equation
\begin{equation*}
\mu_{CZ}(\bar{z})=\Delta_{\psi'}(\bar{z})
\end{equation*}
where $\Delta_{\psi'}(\bar{z})$ is the mean index \cite{SaZ}. This implies that 
\begin{equation*}
\mu_{CZ}(\bar{z}^m)=m\mu_{CZ}(\bar{z})
\end{equation*}
holds for any ${k\in \mathbb{N}}$. For the iteration ${\psi=\psi'^{(2N)}}$ and any capped periodic orbit of ${\psi}$, the same equation ${\mu_{CZ}(\bar{z})=\Delta_{\psi}(\bar{z})}$ holds. Let ${\{y_1,\cdots,y_l\}}$ be all contractible periodic orbits of ${\psi}$. Note that they are ${2N}$-times iterations of ${\{x_1^{(\frac{k}{p_1})},\cdots,x_l^{(\frac{k}{p_l})}\}}$. This means that any capped periodic orbit ${\bar{z}}$ of ${\{y_1,\cdots,y_l\}}$ has mean index ${\Delta_{\psi}(\bar{z})=2Nm}$ (${m\in \mathbb{Z}}$). So we can choose a capping ${\bar{y_i}}$ of ${y_i}$ (${i=1,\cdots,l}$) so that 
\begin{equation*}
\mu_{CZ}(\bar{y_i})=0
\end{equation*}
is satisfied. This implies that Conley-Zehnder index of any capped periodic orbit is divided by ${2N}$ and ${HF_{odd}(\psi)=0}$ holds. If $n$ is an odd integer, this is a contradiction because ${HF_n(\psi)\neq 0}$ holds. So, ${\mathcal{H}_{sn}^{(1)}}$ is empty if ${n}$ is odd.

Next, assume that ${H_{odd}(M:\mathbb{Q})\neq 0}$. Without loss of generality, we assume that $n$ is an even integer. Then the isomorphism ${HF_*(\psi)\cong H_{*-n}(M:\mathbb{Q})\otimes \Lambda_{(M,\omega)}}$ implies that there is at least one capped periodic orbit of ${\psi}$ whose Conley-Zehnder index is odd. This is a contradiction. So ${\mathcal{H}_{sn}^{(1)}}$ is empty if ${H_{odd}(M:\mathbb{Q})\neq 0}$ holds.

Assume that ${N>1}$ holds. Without loss of generality, we assume that ${n}$ is even. Note that ${H_{n+2}(M:\mathbb{Q})\neq 0}$ holds in this case. This implies ${HF_2(\psi)\neq 0}$, but this is impossible because Conley-Zehnder index of any capped periodic periodic orbit can be divided by ${2N}$. So we have proved that ${\mathcal{H}_{sn}^{(1)}}$ is empty in all cases.

Next we fix ${\phi \in \mathcal{H}_{sn}^{(2)}}$. Let ${\{x_1,\cdots,x_l\}}$ be the set of all simple periodic orbits and let $p_1,\cdots,p_l$ be their periods. We divide ${\{x_1,\cdots,x_l\}}$ into hyperbolic periodic orbits and non-hyperbolic periodic orbits. Let ${\{x_1,\cdots,x_{l'}\}}$ be the set of all non-hyperbolic periodic orbits. We perturb ${\phi}$ to ${\widetilde{\phi}}$ so that all differential maps

\begin{equation*}
(d\phi)^{p_i}:T_{x_i(0)}M\longrightarrow T_{x_i(0)}M
\end{equation*}
have ${2n}$ pairwise distinct eigenvalues (see the arguments in the proof of Lemma 7.1.5 in \cite{AD}). The perturbed ${\widetilde{\phi}}$ may not be strongly non-degenerate and perturbed periodic orbits ${\{\widetilde{x_1},\cdots,\widetilde{x_{l'}}\}}$ may not be non-hyperbolic. We prove that at least one ${\widetilde{x_i}}$ is non-hyperbolic. Note that ${\widetilde{\phi}}$ may have simple contractible periodic orbits more than $l$, but we can assume  that the period of "new" periodic orbit is much greater than $2Nk$. So the existence of "new" periodic orbits does not influence our arguments. Assume that all ${\widetilde{x_i}}$ are hyperbolic. As in the proof of ${\mathcal{H}_{sn}^{(1)}=\emptyset}$, we fix ${\psi=\widetilde{\phi}^{(2N\times k)}}$ where $k$ is a common multiple of ${p_1,\cdots,p_l}$. Then each periodic orbit ${\widetilde{x_i}^{(2N\times \frac{k}{p_i})}}$ of ${\psi}$ has a capping ${\bar{z_i}}$ so that ${\mu_{CZ}(\bar{z_i})=0}$. This is a contradiction as in the proof of ${\mathcal{H}_{sn}^{(1)}=\emptyset}$. So at least one of ${\widetilde{x_i}}$ is non-hyperbolic periodic orbit. 

Let ${\widetilde{x_1}}$ be a non-hyperbolic periodic orbit. We can perturb ${\widetilde{\phi}}$ so that ${\widetilde{x_1}}$ satisfies all required conditions in the statement of Birkhoff-Moser fixed point theorem. These arguments imply that we can choose a sequence ${\{\phi_k\}_{k\in \mathbb{N}}}$ and open neighborhoods $W_k$ of ${\phi_k}$ which satisfy the following conditions.

\begin{itemize}
\item ${\phi_k\longrightarrow \phi}$ in ${C^{\infty}}$-topology
\item Any element of ${W_k}$ satisfies all required conditions in the statement of Birkhoff-Moser fixed point theorem and hence possesses infinitely many simple contractible periodic orbits.
\end{itemize}

We define ${V(\phi)}$ and ${V^{(k)}(\phi)}$ (${k\in \mathbb{N}}$) for ${\phi \in \mathcal{H}_{sn}^{(2)}}$ as follows.

\begin{gather*}
V(\phi)=\bigcup_{i=1}^{\infty}W_i   \\
V^{(k)}(\phi)=V(\phi)
\end{gather*}

Next we fix ${\phi \in \mathcal{H}_{sn}^{(3)}}$. There are the following two possibilities.

\begin{enumerate}
\item There is an open neighborhood $U$ of ${\phi}$ (in ${\textrm{Ham}(M,\omega)}$) such that 
\begin{equation*}
U\cap \mathcal{H}_{sn}\subset \mathcal{H}_{sn}^{(3)}
\end{equation*} 
holds.
\item There is no open neighborhood $U$ as above. In other words, we can choose a sequence ${\{\phi_k\}_{k\in \mathbb{N}}\subset  \mathcal{H}_{sn}^{(2)}}$  such that ${\phi_k\rightarrow \phi}$ holds.
\end{enumerate}

In the case of ${(2)}$, we define ${V(\phi)}$ and ${V^{(k)}(\phi)}$ (${k\in \mathbb{N}}$) as follows.

\begin{gather*}
V(\phi)=\bigcup_{k=1}^{\infty}V(\phi_k)  \\
V^{(k)}(\phi)=V(\phi)
\end{gather*}

In the case of ${(1)}$, we define ${V(\phi)}$ and ${V^{(k)}(\phi)}$ (${k\in \mathbb{N}}$) as follows.

\begin{gather*}
V(\phi)=U\cap \mathcal{H}_{sn}  \\
V^{(k)}(\phi)=\{\psi \in U \ | \ \psi,\cdots,\psi^k  \textrm{ \ are\ non-degenerate}\}
\end{gather*}

Then ${V^{(k)}(\phi)\subset U}$ is open dense and ${V(\phi)=\bigcap_{k=1}^{\infty} V^{(k)}(\phi)}$ holds. We can define ${X}$ and ${Y_k}$ as follows.

\begin{gather*}
X=\bigcup_{\phi\in \mathcal{H}_{sn}}V(\phi) \\
Y_k=\bigcup_{\phi\in \mathcal{H}_{sn}}V^{(k)}(\phi)
\end{gather*}

${X}$ and ${\{Y_k\}}$ satisfy the required conditions and we proved Theorem 1.1. 
\begin{flushright}   $\Box$  \end{flushright}

\end{document}